\documentclass[11pt]{article}
\usepackage{amsmath,amsthm,amsfonts,amssymb,latexsym}
\usepackage{amsmath}
\usepackage{color}
\usepackage{a4wide}
\newtheorem{thm}{Theorem}[section]

\newtheorem{lem}[thm]{Lemma}
\newtheorem{prop}[thm]{Proposition}

\numberwithin{equation}{section}
\newtheorem{theorem}{Theorem}[section]

\newtheorem{lemma}[theorem]{Lemma}

\newtheorem{proposition}[theorem]{Proposition}
\newtheorem{corollary}[theorem]{Corollary}

\def\g{\mathfrak g}

\def\a{\mathfrak a}
\def\k{\mathfrak k}

\def\a{\alpha}

\def\l{\lambda}

\def\g{\mathfrak{g}}

\def\k{\mathfrak{k}}

\def\L{{\cal   L}}

\def\HH{{\rm   H}}
\def\fF{{\frak   F}}

\def\cK{{\cal   K}}

\def\half{{\frac{1}{2} }}
\def\cF{{\cal   F}}

\def\L1#1{L^1(#1)}

\def\lef({\left(}
\def\rig){\right)}

\begin{document}

\title{  Cohomology of  $\mathfrak {aff}(1)$ and  $\mathfrak {aff}(1|1)$ acting
on the space of
$n$-ary differential operators on the superspace
$\mathbb{R}^{1|1}$ }

\label{firstpage}

\author{Mabrouk Ben Ammar \and Maha Boujelben \and Amina Jabeur \and Rabeb Sidaoui  \thanks{ Universit\'e de Sfax, Facult\'e des Sciences, D\'epartement de Math\'ematiques, Laboratoire d'Alg\`ebre, G\'eom\'etrie et Th\'eorie Spectrale (AGTS) LR11ES53, BP 802, 3038 Sfax, Tunisie.
E.mail: mabrouk.benammar@fss.rnu.tn }}

\maketitle

\begin{abstract}
 We consider the $\mu$-densities spaces $\mathcal{F}_\mu$ with $\mu\in\mathbb{R}$, we  compute the space $\mathrm{H}^1_\mathrm{diff}(\mathfrak{aff}(1),\mathrm{D}_{\lambda,\mu})$  where $\lambda=(\lambda_1,\dots,\lambda_n)\in\mathbb{R}^n$ and $\mathrm{D}_{\lambda,\mu}$ is the space of $n$-ary differential operators from $\mathcal{F}_{\lambda_1}\otimes\cdots\otimes\mathcal{F}_{\lambda_n}$ to $\mathcal{F}_\mu$. We also compute the super analog space $\mathrm{H}^1_\mathrm{diff}(\mathfrak{aff}(1|1),\mathfrak{D}_{\lambda,\mu})$.
\end{abstract}

\maketitle {\bf Mathematics Subject Classification} 17B56, 17B10, 17B66


\thispagestyle{empty}

\section{Introduction}

The space of weighted densities of weight $\mu$ on
$\mathbb{R}$ (or $\mu$-densities for short), denoted by:
\begin{equation*}{\cal F}_\mu=\left\{ fdx^{\mu}, ~f\in
C^\infty(\mathbb{R})\right\},\quad \mu\in\mathbb{R},
\end{equation*}is the space of
sections of the line bundle $(T^*\mathbb{R})^{\otimes^\mu}.$ The
Lie algebra $\mathrm{Vect}(\mathbb{R})$ of vector fields $X_h=h{d\over
dx}$, where $h\in C^\infty(\mathbb{R})$, acts by the {\it Lie
derivative}. Alternatively, this action can be written as follows:
\begin{equation}\label{Lie1}{X_h}\cdot(fdx^{\mu})=L_{X_h}^\mu(f)dx^{\mu} \quad\text{with}~~L_{X_h}^\mu(f)=hf'+\mu h'f,
\end{equation}
where $f'$, $h'$ are $\frac{df}{dx}$, $\frac{dh}{dx}$. For $\lambda=(\lambda_1,\dots,\lambda_n)\in\mathbb{R}^n$ and $\mu\in\mathbb{R}$ we denote by $\mathrm{D}_{\lambda,\mu}$ the space of $n$-ary differential operators $A$ from $\mathcal{F}_{\lambda_1}\otimes\cdots\otimes\mathcal{F}_{\lambda_n}$ to $\mathcal{F}_\mu$. The Lie algebra $\mathrm{Vect}(\mathbb{R})$ acts
on the space $\mathrm{D}_{\lambda,\mu}$ of these differential
operators by:
\begin{equation}\label{Lieder2}X_h\cdot A=L_{X_h}^\mu\circ
A-A\circ L_{X_h}^{\lambda} \end{equation} where
$L_{X_h}^{\lambda}$ is the Lie derivative on
$\mathcal{F}_{\lambda_1}\otimes\cdots\otimes\mathcal{F}_{\lambda_n}$  defined by the Leibnitz rule.
If we restrict ourselves to the Lie algebra $\mathfrak{ aff}(1)$ which
is isomorphic to the Lie subalgebra of $\mathrm{Vect}(\mathbb{R})$
spanned by
\begin{equation*}\left\{X_1,\,X_x\right\},\end{equation*}
we have a family of infinite dimensional $\frak{ aff}(1)$-modules
still denoted by $\mathrm{D}_{\lambda,\mu}$.

According to Nijenhuis-Richardson \cite{r1}, the space
$\mathrm{H}^1\left(\mathfrak{g};\mathrm{End}(V)\right)$ classifies
the infinitesimal deformations of a $\mathfrak{g}$-module $V$ and
the obstructions to integrability of a given infinitesimal
deformation of $V$ are elements of
$\mathrm{H}^2\left(\mathfrak{g};\mathrm{End}(V)\right)$. While the
spaces $\mathrm{H}^1\left(\mathfrak{g};L(\otimes_s^k V,V)\right)$
appear naturally in the problem of normalization of nonlinear
representations of $\mathfrak{g}$ in $V$. To be more precise, let
$$T:\mathfrak{g}\rightarrow\sum_{k\geq0}L(\otimes_s^k
V,V),~X\mapsto T_X=\sum T_X^k,$$ be a nonlinear representation of
$\mathfrak{g}$ in $V$, that is, $T_{[X,Y]}=[T_X,T_Y]$. In
\cite{abs}, it is proved that the representation $T$ is
normalized if $T_X^k$ is in a supplementary of
$\mathrm{B}^1\left(\mathfrak{g};L(\otimes_s^k V,V)\right)$ in
$\mathrm{Z}^1\left(\mathfrak{g};L(\otimes_s^k V,V)\right)$.

In fact if $A$ is a differential operator on the line, $A$ can be
viewed as an homomorphism from $\mathcal{F}_\lambda$ to
$\mathcal{F}_\mu$. If $A$ is with order $n$, we can define its
symbol as an element in
$\mathcal{S}_\delta^n=\bigoplus_{j=0}^n\mathcal{F}_{\delta-j}$ for
$\delta=\mu-\lambda$. If $n$ goes to $+\infty$, the space
$\mathcal{S}_\delta=\bigoplus_{j\geq0}\mathcal{F}_{\delta-j}$
appears as the space of symbols for all differential operators. The
space
$\mathrm{H}^1\left(\mathfrak{g};L(\otimes_s^2
\mathcal{S}_\delta,\mathcal{S}_\delta)\right)$ can be decomposed as
a sum of spaces $\mathrm{H}^1\left(\mathfrak{g},
\mathrm{D}_{\lambda,\mu}\right)$ with $\lambda=(\lambda_1,\lambda_2)\in\mathbb{R}^2$. Thus, the computation of the
spaces $\mathrm{H}^1\left(\mathfrak{g},
\mathrm{D}_{\lambda,\mu}\right)$ is the first step to
normalize any nonlinear representation of $\mathfrak{g}$ in
$\mathcal{S}_\delta$.

For $\lambda\in\mathbb{R}$ the spaces $\mathrm{H}^1_\mathrm{diff}\left(\mathfrak{sl}(2),
\mathrm{D}_{\lambda,\mu}\right)$ are computed by Gargoubi \cite{g} and Lecomte \cite{lec} and the spaces $\mathrm{H}^1_\mathrm{diff}\left(\mathfrak{osp}(1|2),
\mathfrak{D}_{\lambda,\mu}\right)$ are computed by Basdouri and Ben Ammar \cite{bab}, where $\mathrm{H}^1_\mathrm{diff}$ denotes the differential cohomology;
that is, only cochains given by differential operators are
considered. For $\lambda\in \mathbb{R}^2$ the spaces $\mathrm{H}^1_\mathrm{diff}\left(\mathfrak{sl}(2),
\mathrm{D}_{\lambda,\mu}\right)$ are computed by Bouarroudj \cite{b} and the spaces $\mathrm{H}^1_\mathrm{diff}\left(\mathfrak{osp}(1|2),
\mathfrak{D}_{\lambda,\mu}\right)$ are computed by Ben Ammar et al \cite{bjs}, while we are interested in this paper in the spaces $\mathrm{H}^1_\mathrm{diff}\left(\mathfrak{aff}(1),
\mathrm{D}_{\lambda,\mu}\right)$ and $\mathrm{H}^1_\mathrm{diff}\left(\mathfrak{aff}(1|1),
\mathfrak{D}_{\lambda,\mu}\right)$ where $\lambda\in\mathbb{R}^n$.
\section{Definitions and Notations}
\subsection{The Lie superalgebra of contact vector fields on $\mathbb{R}^{1|1}$}
We define the superspace $\mathbb{R}^{1|1}$ in terms of its
superalgebra of functions, denoted by $C^\infty(\mathbb{R}^{1|1})$
and consisting of elements of the form:
\begin{equation*}
F(x,\theta)=f_0(x)+f_1(x)\theta ,
\end{equation*}
where $x$ is the even variable,  $\theta$ is  the odd variable
($\theta^2=0$) and $f_0(x),\,f_1(x)\in C^\infty(\mathbb{R})$. Even elements in $C^\infty(\mathbb{R}^{1|1})$ are the functions
$F(x,\theta)=f_0(x)$, the functions $F(x,\theta)=\theta f_1(x)$ are odd elements. We
consider the contact bracket on $C^\infty(\mathbb{R}^{1|1})$
defined on  $C^\infty(\mathbb{R}^{1|1})$ by:
\begin{equation*}\begin{array}{l}
\{F,G\}=FG'-F'G+\frac{1}{2}D(F)
\overline{D}(G),\end{array}
\end{equation*}where $D=\frac{\partial}{\partial
{\theta}}+\theta\frac{\partial}{\partial x}$ and
$\overline{D}=\frac{\partial}{\partial
{\theta}}-\theta\frac{\partial}{\partial x}$. The superspace $\mathbb{R}^{1|1}$ is equipped with
the standard contact structure given by the following $1$-form:
\begin{equation*}
\a=dx+\theta d\theta.
\end{equation*}
Let $\mathrm{Vect}(\mathbb{R}^{1|1})$ be the superspace of
 vector fields on $\mathbb{R}^{1|1}$:
\begin{equation*}\mathrm{Vect}(\mathbb{R}^{1|1})=\left\{F_0\partial_x
+  F_1\partial_\theta \mid ~F_i\in C^\infty(
\mathbb{R}^{1|1})\right\},\end{equation*} where $\partial_\theta$
stands for $\frac{\partial}{\partial\theta}$ and $\partial_x$
stands for $\frac{\partial}{\partial x} $, and consider the
superspace $\mathcal{K}(1)$ of contact vector fields on
$\mathbb{R}^{1|1}$. That is, $\mathcal{K}(1)$ is the superspace
of vector fields on $\mathbb{R}^{1|1}$ preserving the distribution
$\langle\overline{D}\rangle$:
$$
\mathcal{K}(1)=\big\{X\in\mathrm{Vect}(\mathbb{R}^{1|1})~|~[X,\,\overline{D}]=
F_X\overline{D}\quad\hbox{for some}~F_X\in
C^\infty(\mathbb{R}^{1|1})\big\}. $$
Any contact vector field
on $\mathbb{R}^{1|1}$ has the following explicit form:
\begin{equation*}\begin{array}{l}
X_F=F\partial_x+\half D(F)\overline{D},\;\text{ where }\, F\in
C^\infty\left(\mathbb{R}^{1|1}\right).\end{array}
\end{equation*}
The bracket on $\mathcal{K}(1)$ is given by
\begin{equation*}
[X_{F},\,X_{G}]=X_{\{F,\,G\}}.
\end{equation*}
Thus, the map  $F\mapsto X_F$ is a Lie superalgebra isomorphism from $\cK(1)$  to
$C^\infty\left(\mathbb{R}^{1|1}\right)$.
Of course, the set $\left\{X_f=f\partial_x+\frac{1}{2}f'\theta\partial\theta,\, f\in C^\infty(\mathbb{R})\right\}$
is a subalgebra of $\mathcal{K}(1)$ isomorphic to $\mathrm{Vect}(\mathbb{R})$.
\subsection{The subalgebra $\mathfrak{aff}(1|1)$}

The Lie algebra $\mathfrak {aff}(1)$ is realized as subalgebra of the
Lie algebra  $\mathrm{Vect}(\mathbb{R})$:
\begin{equation*}\frak {aff}(1)=\text{Span}(X_1,\,X_{x}).\end{equation*}
Similarly, we now consider the Lie superalgebra $\mathfrak {aff}(1|1)$ as
a subalgebra of $\mathcal{K}(1)$:
\begin{equation*}
\mathfrak{aff}(1|1)=\text{Span}(X_1,\,X_{x},\,
X_{\theta}).
\end{equation*}
The space of even elements is isomorphic to $\mathfrak{aff}(1)$, while the space of odd elements is one dimensional:
\begin{equation*}
(\mathfrak{aff}(1|1))_{\bar1}=\text{Span}(X_{\theta}).
\end{equation*}
The commutation relations are
\begin{equation*}
\begin{array}{llll}
&[X_x,X_\theta]~=-\frac{1}{2}
X_\theta,~~&[X_x,X_1]=-X_1,
&[X_{\theta},X_{\theta}]~=\frac{1}{2}X_{1}.
\end{array}
\end{equation*}

\subsection{The space of weighted densities on $\mathbb{R}^{1|1}$}

We have analogous definition of weighted
densities in super setting (see \cite{ab}) with $dx$ replaced by
$\alpha$. The elements of these spaces are indeed (weighted)
densities since all spaces of generalized tensor fields have just
one parameter relative $\mathcal{K}(1)$
--- the value of $X_x$ on the lowest weight vector (the one
annihilated by $X_\theta$). From this point of view the volume
element (roughly speaking, $\lq\lq
dx\frac{\partial}{\partial\theta}"$) is indistinguishable from
$\alpha^{\frac{1}{2}}.$ We
denote by $\mathfrak{F}_{\mu}$ the space of all
 weighted densities on $\mathbb{R}^{1|1}$ of weight $\mu\in\mathbb{R}$:
\begin{equation*}
\mathfrak{F}_\mu=\left\{F(x,\theta)\alpha^\mu~~|~~F(x,\theta)
\in C^\infty(\mathbb{R}^{1|1})\right\}.
\end{equation*}
As a vector space, $\mathfrak{F}_\mu$ is isomorphic to
$C^\infty(\mathbb{R}^{1|1})$, but the Lie derivative of the
density $G\alpha^\mu$ along the vector field $X_F$ in
$\mathcal{K}(1)$ is now:
\begin{equation*}
\label{superaction}
\mathfrak{L}_{X_F}(G\alpha^\mu)=\mathfrak{L}^{\mu}_{X_F}(G)\alpha^\mu,
\quad\text{with}~~\mathfrak{L}^{\mu}_{X_F}(G)=\mathfrak{L}_{X_F}(G)+
\mu F'G.
\end{equation*}
Especially, if  $f\in C^\infty(\mathbb{R})$ and
$G(x,\theta)=g_0(x)+g_1(x)\theta$, then we easily check that
\begin{equation}
\label{deriv}
\mathfrak{L}^{\mu}_{X_f}(G)=L^{\mu}_{X_f}(g_0)+
\left(L^{\mu+\frac{1}{2}}_{X_f}(g_1)\right)\theta.
\end{equation}
Of course, for all $\mu$, $\mathfrak{F}_\mu$ is a
$\mathcal{K}(1)$-module:
\begin{equation*}
[\mathfrak{L}^{\mu}_{X_F},\mathfrak{L}^{\mu}_{X_G}] =
\mathfrak{L}^{\mu}_{[{X_F},\,X_G]}.
\end{equation*}

\subsection{Differential operators on weighted densities}

A differential operator on $\mathbb{R}^{1|1}$ is an operator on
$C^\infty(\mathbb{R}^{1|1})$ of the following form:
$$
A= \sum_{i=0}^\ell{a}_i (x, \theta)\partial_x^i +
\sum_{i=0}^\ell{b}_i (x,
\theta)\partial_x^i\partial_\theta.
$$
In \cite{gmo}, it is proved that any local operator $A$ on
$\mathbb{R}^{1|1}$ is in fact a differential operator.

Of course, any differential operator defines a linear mapping from
$\mathfrak{F}_\lambda$ to $\mathfrak{F}_\mu$ for any $\lambda$,
$\mu\in\mathbb{R}$: $F\alpha^\lambda\mapsto A(F)\alpha^\mu$. Similarly, if $\lambda=(\lambda_1,\dots,\lambda_n)\in\mathbb{R}^n$ and $\mu\in\mathbb{R}$, we consider  a family of $\mathcal{K}(1)$ and $\mathfrak{aff}(1|1)$
modules denoted $\mathfrak{D}_{\lambda,\mu}$, for the natural
action:
\begin{equation*}\label{d-action}
{X_F}\cdot A=\mathfrak{L}^\mu_{X_F}\circ A-
(-1)^{AF}A\circ \mathfrak{L}^\lambda_{X_F}
\end{equation*}
where ${\frak
L}_{X_F}^{\lambda}$ is the Lie derivative on ${\frak
F}_{\lambda_1}\otimes\cdots\otimes\frak{F}_{\lambda_n}$ defined by the Leibnitz rule.

\section{The space $\HH^1_\mathrm{diff}(\frak {aff}(1|1),\frak{D}_{\lambda,\mu})$}

\subsection{Cohomology }
We will compute the first cohomology space of $\mathfrak {aff}(1|1)$
with coefficients in $\mathfrak{D}_{\lambda,\mu}$ where $\lambda\in\mathbb{R}^n$ and $\mu\in\mathbb{R}$. Let us first
recall some fundamental concepts from cohomology theory~(see, e.g.,
\cite{Fu}). Let $\mathfrak{g}=\frak{g}_{\bar 0}\oplus \mathfrak{g}_{\bar 1}$
be a Lie superalgebra acting on a superspace $V=V_{\bar 0}\oplus
V_{\bar 1}$ and let $\mathfrak{h}$ be a subalgebra of
$\mathfrak{g}$. (If $\frak{h}$ is omitted it assumed to be $\{0\}$).
The space of $\frak h$-relative $n$-cochains of $\mathfrak{g}$ with
values in $V$ is the $\mathfrak{g}$-module
\begin{equation*}
C^n(\mathfrak{g},\mathfrak{h}; V ) := \mathrm{Hom}_{\mathfrak
h}(\Lambda^n(\mathfrak{g}/\mathfrak{h});V).
\end{equation*}
The {\it coboundary operator} $ \partial^n: C^n(\mathfrak{g},\mathfrak{h}; V
)\longrightarrow C^{n+1}(\mathfrak{g},\mathfrak{h}; V )$ is a
$\mathfrak{g}$-map satisfying $\partial^n\circ\partial^{n-1}=0$. The
kernel of $\partial^n$, denoted $Z^n(\mathfrak{g},\frak{h};V)$, is
the space of $\mathfrak{h}$-relative $n$-{\it cocycles}, among them,
the elements in the range of $\partial^{n-1}$ are called $\mathfrak{h}$-relative $n$-{\it coboundaries}. We denote
$B^n(\mathfrak{g},\mathfrak{h};V)$ the space of $n$-coboundaries.

By definition, the $n^{th}$ $\mathfrak h$-relative  cohomology space is
the quotient space
\begin{equation*}
\mathrm{H}^n
(\mathfrak{g},\frak{h};V)=Z^n(\mathfrak{g},\frak{h};V)/B^n(\mathfrak{g},\frak{h};V).
\end{equation*}
We will only need the formula of $\partial^n$ (which will be simply
denoted $\partial$) in degrees 0 and 1. For $v \in
C^0(\mathfrak{g},\,\mathfrak{h}; V) =V^{\mathfrak h}$, $$\partial v(g) : =
(-1)^{gv}g\cdot v,$$ where $V^{\mathfrak h}$ is the subspace of $\mathfrak{h}$-invariant elements of $V$. For  $ \Omega\in C^1(\mathfrak{g}, \mathfrak{h};V )$ and $ g,\,h\in\mathfrak{g}$,
\begin{equation*}\partial(\Omega)(g,\,h):=
(-1)^{g\Omega}g\cdot
\Omega(h)-(-1)^{h(g+\Omega)}h\cdot
\Omega(g)-\Upsilon([g,~h]).
\end{equation*}
\begin{proposition}\label{pr}\begin{itemize}
                           \item [1)] Let $\Omega\in \mathrm{Z}^1(\mathfrak{g};V)$ and $a\in\mathfrak{g}$. If $\Omega(a)=0$ then
$a\cdot \Omega=0$. Thus, if $\Omega|_\mathfrak{h}=0$ then $\Omega$ is $\mathfrak{h}$-invariant: $\Omega\in \mathrm{Z}^1(\mathfrak{g},\mathfrak{h};V)$. Moreover, $\mathrm{H}^1(\mathfrak{g},
\mathfrak{h};V)\subset\mathrm{H}^1(\mathfrak{g};V)$.
                           \item [2)] If $\mathrm{H}^1(\mathfrak{g},
\mathfrak{h};V)=0$ then $\Omega\in\mathrm{B}^1(\mathfrak{g};V)$ if and only if $\Omega|_\mathfrak{h}\in\mathrm{B}^1(\mathfrak{h};V)$. Thus, in this case, the space $\mathrm{H}^1(\mathfrak{g};V)$ is characterized by the space $\mathrm{H}^1(\mathfrak{h};V)$.
                         \end{itemize}
 \end{proposition}
\begin{proofname}. 1) Indeed, for any $x\in\mathfrak{g}$, we have
$$
(-1)^{a\Omega}a\cdot \Omega(x)-(-1)^{x(a+\Omega)}x\cdot \Omega(a)-\Omega([a,~x])=0\,\text{and}\,(a\cdot \Omega)(x)=a\cdot \Omega(x)-(-1)^{a\Omega}\Omega([a,~x]).
$$
Moreover, for $\Omega\in\mathrm{Z}^1(\mathfrak{g},
\mathfrak{h};V)$ we denote by $\overline{\Omega}$ (respectively $\widetilde{\Omega}$) the class of $\Omega$ up to $\mathrm{B}^1(\mathfrak{g},
\mathfrak{h};V)$ (respectively $\mathrm{B}^1(\mathfrak{g};V)$). We easily check that the map $\overline{\Omega}\mapsto\widetilde{\Omega}$, from $\mathrm{H}^1(\mathfrak{g},
\mathfrak{h};V)$ to $\mathrm{H}^1(\mathfrak{g};V)$, is injective : if $\Omega=\partial v$ with $v\in V$ then $v\in V^\mathfrak{h}$ since $\Omega|_\mathfrak{h}=0$.

2) Obviously if $\Omega\in\mathrm{B}^1(\mathfrak{g};V)$ then $\Omega|_\mathfrak{h}\in\mathrm{B}^1(\mathfrak{h};V)$. Inversely, if $\Omega|_\mathfrak{h}\in\mathrm{B}^1(\mathfrak{h};V)$ then, up to a coboundary, $\Omega$ vanishes on $\mathfrak{h}$, therefore,  $\Omega\in\mathrm{Z}^1(\mathfrak{g},\mathfrak{h};V)$. But $\mathrm{H}^1(\mathfrak{g},\mathfrak{h};V)=\{0\}\subset\mathrm{H}^1(\mathfrak{g};V)$, so, $\Omega\in\mathrm{B}^1(\mathfrak{g};V)$.\hfill$\Box$\end{proofname}

\subsection{The space  $\HH^1_\mathrm{diff}(\frak
{aff}(1),\mathrm{D}_{\lambda,\mu})$}

 Let $\mu\in\mathbb{R}$, $\alpha=(\alpha_1\dots,\alpha_n)\in\mathbb{N}^n$ and $\lambda=(\lambda_1\dots,\lambda_n)\in\mathbb{R}^n$, we consider $\delta=\mu-\sum_{i=1}^n\lambda_i$ and $|\alpha|=\sum\alpha_i.$
 For $F=f_1\otimes\cdots\otimes f_n\in\mathcal{F}_{\lambda_1}\otimes\cdots \otimes\mathcal{F}_{\lambda_n} $, we define $$F^{(\alpha)}:=f_1^{(\alpha_1)}\dots f_n^{(\alpha_n)}.$$
  Recall that the space $\mathcal{F}_{\lambda_1}\otimes\cdots \otimes\mathcal{F}_{\lambda_n} $ is an $\mathfrak{aff}(1)$-module:
$$
 X_h\cdot F:=L_{X_h}^{\lambda}(F)=\sum_{i=1}^nf_1\otimes\cdots\otimes L_{X_h}^{\lambda_i}(f_i)\otimes\cdots\otimes f_n.
$$

The following lemma gives the general form of any 1-cocycle.
\begin{lem}\label{lem1} Up to a coboundary, any 1-cocycle
$\Omega\in\mathrm{Z}^1_\mathrm{diff}(\frak{aff}(1),\mathrm{D}_{\lambda,\mu})$
can be expressed as follows. For all $F=f_1\otimes\cdots\otimes f_n\in \mathcal{F}_{\lambda_1}\otimes\cdots \otimes\mathcal{F}_{\lambda_n} $ and for all $X_h\in\frak{aff}(1)$:
\begin{equation}\label{cocycle}\Omega(X_h,F) =\sum_\alpha B_\alpha h' F^{(\alpha)} ,\end{equation}
where  the  $B_\alpha$ are constants.

\end{lem}

\begin{proofname}. Any 1-cocycle on $\frak{aff}(1)$ should retains the following
general form:
$$\Omega(X_h,F) = \sum_\alpha N_\alpha h F^{(\alpha)}+\sum_\alpha M_\alpha h' F^{(\alpha)} ,$$
 where $N_\alpha$ and $M_\alpha$ are, a priori, functions. First, we prove that
the terms in $h$ can be annihilated by adding a coboundary.
Let $b : \mathcal{F}_{\lambda_1}\otimes\cdots \otimes\mathcal{F}_{\lambda_n}\rightarrow
\mathcal{F_\mu}$ be a $n$-ary differential operator defined by
$$
b(F) =
\sum_\alpha D_\alpha F^{(\alpha)} ,
$$
We have
\begin{equation}\label{delta}\begin{array}{lll}
\partial b(X_h,F) &=&h(b(F))'+\mu h'b(F)-b(X_h\cdot F)\\[10pt]
&=& \displaystyle\sum_\alpha D'_\alpha h F^{(\alpha)}+\displaystyle\sum_\alpha (\delta-|\alpha|)D_\alpha h' F^{(\alpha)}
\end{array}\end{equation}

Thus, if $D'_\alpha=N_\alpha$ then $\Omega-\partial b$ does not contain terms in h. So, we can replace $\Omega$ by
$\Omega-\partial b$.
 That is, up to a
coboundary, any 1-cocycle on $\frak{aff}(1)$ can
be expressed as follows:
$$\Omega(X_h,F) = \sum_\alpha B_\alpha h' F^{(\alpha)}.$$
Now, consider the 1-cocycle condition:
$$
\Omega([X_{h_1}, X_{h_2}],F)- X_{h_1}\cdot \Omega(X_{h_2}, F)
+X_{h_2}\cdot \Omega(X_{h_1}, F) = 0,
$$
where $X_{h_1},\,X_{h_2} \in\frak{aff}(1)$. That is,
\begin{equation*}\begin{array}{l} \displaystyle\sum_\alpha B'_\alpha ( h_1 h'_2-h'_1h_2) F^{(\alpha)}=0.
\end{array}\end{equation*}
 So, for all $\alpha$, we have $
B'_\alpha=0$
\hfill$\Box$\end{proofname}

\begin{thm}\label{cor} \begin{itemize}
                         \item [1)] If $\delta\notin\mathbb{N}$ then ${\rm H}^1_\mathrm{diff}(\frak{aff}(1),\mathrm{D}_{\lambda,\mu})=0$.
                         \item [2)] If $\delta\in\mathbb{N}$ then, up to a coboundary,  any 1-cocycle
$c\in\mathrm{Z}^1_\mathrm{diff}(\frak{aff}(1),\mathrm{D}_{\lambda,\mu})$
can be expressed as follows. For all $F=f_1\otimes\cdots\otimes f_n\in \mathcal{F}_{\lambda_1}\otimes\cdots \otimes\mathcal{F}_{\lambda_n} $ and for all $X_h\in\frak{aff}(1)$:
\begin{equation}\label{cocyc}\Omega(X_h,F) =\sum_{|\alpha|=\delta} B_\alpha h' F^{(\alpha)},\end{equation}

\end{itemize}

\end{thm}
\begin{proofname}. 1) Indeed, according to Lemma \ref{lem1}, we can easily show  the 1-cocycle $\Omega$ defined by (\ref{cocycle}) is nothing but the operator
$\partial b$ where
$$
b(F) =
\sum_{\alpha} {B_\alpha\over \delta-|\alpha|} F^{(\alpha)} ,
$$

2)
Consider the 1-cocycle $\Omega$ defined by (\ref{cocycle}) and consider the operator
$\partial b$ where
$$
b(F) =
\sum_{|\alpha|\neq \delta} {B_\alpha\over \delta-|\alpha|} F^{(\alpha)}.
$$
 We easily show that
$$
(\Omega-\partial b)(X_h,F)=\sum_{|\alpha|=\delta} B_\alpha h' F^{(\alpha)}.
$$

\hfill$\Box$\end{proofname}
\begin{thm}  If $\delta=k\in\mathbb{N}$ then
 $$\mathrm{dim}\mathrm{H}^1(\mathfrak{aff}(1);\mathrm{D}_{\lambda,\mu})= \begin{pmatrix}n+k-1\\k\end{pmatrix}.$$
 \end{thm}
 \begin{proofname}.
 According to Theorem \ref{cor} and to the formula \eqref{delta}, the classes of 1-cocycles $\Omega^\alpha$ defined by  $\Omega^\alpha(X_h,F)= h' F^{(\alpha)}$, where $|\alpha|=k$, constitute a basis of $\mathrm{H}^1(\mathfrak{aff}(1);\mathrm{D}_{\lambda,\mu})$. Thus, $\mathrm{dim}\mathrm{H}^1(\mathfrak{aff}(1);\mathrm{D}_{\lambda,\mu})$ is the cardinal of the set $\left\{\alpha\in \mathbb{N}^n,~|\alpha|=k\right\}$.
 \end{proofname}

 \section{Relationship between $\HH^1_\mathrm{diff}(\frak
{aff}(1|1),\frak{D}_{\l,\mu})$ and $\HH^1_\mathrm{diff}(\frak
{aff}(1),\mathrm{D}_{\lambda,\mu})$}
We need to present here some results illustrating the analogy
between the cohomology spaces in super and classical settings.
We consider $$\begin{array}{l}\mathcal{E}=\left\{\varepsilon=(\varepsilon_1,\dots,\varepsilon_{n}),\quad \varepsilon_i=0,\,\half,\quad\text{and}\quad|\varepsilon|\in\mathbb{N}\right\},\end{array}$$
$$\begin{array}{l}\mathcal{O}=\left\{\varepsilon=(\varepsilon_1,\dots,\varepsilon_{n}), 
\quad \varepsilon_i=0,\,\half,\quad\text{and}\quad|\varepsilon|\notin\mathbb{N}\right\}.\end{array}$$
\begin{proposition}\label{prop1}~

\noindent 1) As a $\frak {aff}(1)$-module, we have
$$\fF_{\mu}\simeq\cF_\mu \oplus \Pi(\cF_{\mu+\half})\quad\text{ and }\quad \frak
{aff}(1|1)\simeq\frak {aff}(1)\oplus \Pi(\mathbb{R}dx^{-\frac{1}{2}}),$$ where
$\mu\in\mathbb{R}$ and $\Pi$ is the
change of parity.

\noindent 2)  For $\lambda\in\mathbb{R}^n$ and $\mu\in\mathbb{R}$, as a $\frak {aff}(1)$-module, we have, for
the homogeneous relative parity components:
\begin{equation}\label{deven}(\frak{D}_{\lambda,\mu})_{\bar0}
\simeq\bigoplus_{\varepsilon\in\mathcal{E}}\mathrm{D}_{(\lambda+\varepsilon,\mu)}\oplus
\bigoplus_{\varepsilon\in\mathcal{O}}\mathrm{D}_{(\lambda+\varepsilon,\mu+\half)}\quad\text{and}\quad
(\frak{D}_{\lambda,\mu})_{\bar1}\simeq\Pi\left(\bigoplus_{\varepsilon\in\mathcal{O}}\mathrm{D}_{(\lambda+\varepsilon,\mu)}\oplus
\bigoplus_{\varepsilon\in\mathcal{E}}\mathrm{D}_{(\lambda+\varepsilon,\mu+\half)}\right).
\end{equation}
\end{proposition}
\begin{proofname}.  The first statement is immediately deduced
 from (\ref{deriv}) and from the fact that $\mathcal{K}(1)\simeq\mathfrak{F}_{-\half}$. The second statement can be deduced from the first one (see \cite{bjs}).

\hfill$\Box$\end{proofname}

Now, in order to describe $\HH^1(\frak
{aff}(1|1),\frak{D}_{\lambda,\mu})$, we need first to describe the $\mathfrak{aff}(1)$-relative cohomology
space $\HH^1_\mathrm{diff}(\mathfrak{aff}(1|1),\frak
{aff}(1);\mathfrak{D}_{\lambda,\mu})$. So, we shall need the following
description of some $\mathfrak{aff}(1)$-invariant
 mappings.
\begin{lem}
\label{inva}Let
\begin{equation}\begin{array}{llll}\label{invar}A:\mathbb{R}dx^{-\half}\otimes\cF_{\lambda_1}\otimes\cdots\otimes\cF_{\lambda_n}&\longrightarrow&\cF_\mu,\\
(adx^{-\half},f_1dx^{\l_1},\dots,f_ndx^{\l_n})&\longmapsto&
A(a,f_1,\dots,f_n)(dx)^{\mu}\end{array}\end{equation} be an $n+1$-ary differential
operator and let $\delta=\mu-\sum_i\l_i$. If $A$ is a nontrivial  $\frak{aff}(1)$-invariant operator then
\begin{equation*}\begin{array}{ll}\delta+\half\in\mathbb{N}.\end{array}\end{equation*}  The corresponding
operator $A$ is  of the form:
\begin{equation*}
A(a,f_1,\dots,f_n)=\sum_{|\alpha|=\delta+\half}c_{\alpha}F^{(\alpha)},
\end{equation*}
where the $c_{\alpha}$ are constants.
\end{lem}
\begin{proofname}. The invariance with respect the vector field $X_1=\partial_x$ yields that $A$ must be expressed with constant coefficients. Thus, the operator $A$ can be expressed as follows: \begin{equation*}
A(a,f_1,\dots,f_n)=\sum_{\alpha}c_{\alpha}F^{(\alpha)},
\end{equation*} where the $c_{\alpha}$ are constants, $\alpha=(\alpha_1,\dots,\alpha_n)\in\mathbb{N}^n$ and $F^{(\alpha)}=f_1^{(\alpha_1)}\cdots f_n^{(\alpha_n)}$.
Note that the operator $A$ can be viewed as element of $D_{\lambda,\mu+\half}$.

The invariance property of $A$ with respect the vector fields $X_x$ reads:
\begin{equation}\label{inv}\begin{array}{llll} 0&=&x(A(f_1,\dots,f_n))'+(\mu+\half ) A(f_1,\dots,f_n)-
A(xf_1'+\lambda_1 f_1,f_2,\dots, f_n)\\&~&-A(f_1,xf_2'+\lambda_2 f_2,f_3,\dots, f_n)-\cdots-A(f_1,f_2,\dots, xf_n'+\lambda_n f_n).\end{array}\end{equation}
Consider any non vanishing coefficient $c_{\alpha}$ and consider terms in
$F^{(\alpha)}$ in (\ref{inv}), we get
\begin{equation*}\begin{array}{ll}\delta+\half=|\alpha|.\end{array}\end{equation*}
\hfill$\Box$ \end{proofname}

\begin{prop}\label{sa} The $\mathfrak{aff}(1)$-relative cohomology
spaces $\HH^1_\mathrm{diff}(\mathfrak{aff}(1|1),\frak
{aff}(1);\mathfrak{D}_{\lambda,\mu})$ are all trivial. That is,
any 1-cocycle $\Omega$ is a coboundary over $\frak{aff}(1|1)$ if
and only if its restriction  to $\frak{aff}(1)$ is a coboundary over
$\frak{aff}(1)$.
\end{prop}
\begin{proofname}. First, it is well known that the space
$\HH^1_\mathrm{diff}(\mathfrak{aff}(1|1),\frak
{aff}(1);\mathfrak{D}_{\lambda,\mu})$ is nothing but the space of
cohomology classes of 1-cocycles vanishing on $\frak{aff}(1)$ \cite{bjs}.

Let $\Omega$ be a 1-cocycle vanishing on $\frak{aff}(1)$, then, by the 1-cocycle condition, we have:
\begin{align}
\label{sltr1} &{X_1}\cdot\Omega(X_{\theta})-\Omega([X_1,X_{\theta}])=0, \\[2pt]
\label{sltr2} &{X_x}\cdot\Omega(X_{\theta})-\Omega([X_x,X_{\theta}])=0, \\[2pt]
\label{sltr3}&{X_{\theta}}\cdot\Omega(X_{\theta})=0.
 \end{align}
 The equations \eqref{sltr1} and \eqref{sltr2} express the
$\frak{aff}(1)$-invariance property of the map $\Omega$. According to the decomposition \eqref{deven} the map $\Omega(X_\theta)$ is decomposed
into some $\frak{aff}(1)$-invariant elements $\Omega$ of $\displaystyle\bigoplus_{\varepsilon\in\mathcal{E}}\mathrm{D}_{(\lambda+\varepsilon,\mu)}\oplus
\bigoplus_{\varepsilon\in\mathcal{O}}\mathrm{D}_{(\lambda+\varepsilon,\mu+\half)}$ or $\displaystyle\Pi\left(\bigoplus_{\varepsilon\in\mathcal{O}}\mathrm{D}_{(\lambda+\varepsilon,\mu)}\oplus
\bigoplus_{\varepsilon\in\mathcal{E}}\mathrm{D}_{(\lambda+\varepsilon,\mu+\half)}\right)$, according to that $\Omega$ is odd or even.  Therefore, the
expressions of these maps are given by Lemma \ref{inva}, in fact,
the change of parity functor $\Pi$ commutes with the the
$\mathfrak{aff}(1)$-action. 

So, if $\Omega$ is an even 1-cocycle then we must have $\delta\in\mathbb{N}$ and if $\Omega$ is an odd 1-cocycle then we must have
$\delta+\half\in\mathbb{N}$. Otherwise, the
operator $\Omega$ is identically the zero map.

If $\Omega$ is an even 1-cocycle then $\Omega(X_\theta)$ is odd, so
$$
\Omega(X_\theta)=\sum_{\varepsilon\in\mathcal{O},|\alpha|=\delta_\varepsilon+\frac{1}{2}}C^{\varepsilon}_\alpha \Omega_{\varepsilon}^\alpha+
\theta\sum_{\varepsilon\in\mathcal{E},|\alpha|=\delta_\varepsilon+1}D^{\varepsilon}_\alpha \Omega_{\varepsilon}^\alpha,
$$
where, 
if $F_i=f_i+\theta g_i$, then
$$
\Omega_\varepsilon^\alpha(F_1,\dots,F_n)= H_\varepsilon^{(\alpha)}
$$
with
$$
 H_\varepsilon^{(\alpha)}=h_1^{(\alpha_1)}\cdots h_n^{(\alpha_n)} \quad\text{where}\quad h_i=\left\{\begin{array}{ll}f_i\quad\text{if}\quad\varepsilon_i=0\\g_i\quad\text{if}\quad\varepsilon_i=\half\end{array}\right.
 $$
\begin{lemma}\label{lem2}
For any $X_h\in\mathfrak{aff}(1)$ we have
$$
\partial\theta\Omega_\varepsilon^\alpha(X_h)=(\delta_\varepsilon+\frac{1}{2}-|\alpha|)\theta\Omega_\varepsilon^\alpha,
$$
where $\Omega_\varepsilon^\alpha$ is viewed as element of $\mathfrak{D}_{\lambda,\mu}$.
\end{lemma}
\begin{proofname}.
By a straightforward computation. We can also see $\Omega_\varepsilon^\alpha$  as element of $\mathrm{D}_{\lambda+\varepsilon,\mu+\frac{1}{2}}$, so, we deduce the result from \eqref{delta}. 
\end{proofname}

Now, since $X_\theta=\frac{1}{2}D$, then by the equation \eqref{sltr3} we have
\begin{equation}\label{theta}\begin{array}{llll} D(\Omega(X_\theta,F_1,\dots,F_n))=-\sum_{i=1}^n\xi_i
\Omega(X_\theta,F_1,\dots,D(F_i),\dots, F_n)\end{array}\end{equation}
where $$\xi_i=(-1)^{\sum_{j=1}^{i-1}F_j}.$$

For $\varepsilon\in\mathcal{E}\cup\mathcal{O}$ and $i\in\{1,\dots,n\}$ we consider the element $\varepsilon^i=\eta\in\mathcal{E}\cup\mathcal{O}$   defined by
$$
\eta_i\neq\varepsilon_i\quad\text{and}\quad\eta_j=\varepsilon_j\quad\text{for}\quad j\neq i
$$
and we define $\xi^i_\varepsilon$  by
$$
\xi^i_\varepsilon=(-1)^{2\sum_{j<i}\varepsilon_j}.
$$
Note that \begin{equation}\label{ei}\xi^i_\varepsilon=\xi^i_{\varepsilon^i}.\end{equation}

For $\alpha\in\mathbb{N}^n$, we define
 $$\alpha^i=(\alpha_1,\dots,\alpha_i+1,\dots, \alpha_n)\quad{ and }\quad \alpha^{\overline{i}}=(\alpha_1,\dots,\alpha_i-1,\dots, \alpha_n).$$
Note that \begin{equation}\label{varepsilon}(\varepsilon^i)^i=\varepsilon,\quad(\alpha^{\overline{i}})^i=(\alpha^{i})^{\overline{i}}=
\alpha\quad\text{and}\quad \varepsilon\in\mathcal{E}\Leftrightarrow\varepsilon^i\in\mathcal{O}.\end{equation}

If $F=f+\theta g$ then $D(F)=g+\theta f'$. Thus,
$$
\Omega_\varepsilon^\alpha(F_1,\dots,D(F_i),\dots,F_n)=\Omega_{\varepsilon^i}^{\alpha_\varepsilon^i}(F_1,\dots,\dots,F_n).
 $$
 where
 $$
 \alpha_\varepsilon^i=\left\{\begin{array}{ll}\alpha\quad\text{if}\quad\varepsilon_i=0\\
 \alpha^i\quad\text{if}\quad\varepsilon_i=\half\end{array}\right..
 $$
 Define $\alpha_\varepsilon^{\overline{i}}$ so that
 \begin{equation}\label{alpha}(\alpha_\varepsilon^{\overline{i}})_{\varepsilon^i}^i=\alpha,\end{equation}
  that is
  $$
 \alpha_\varepsilon^{\overline{i}}=\left\{\begin{array}{ll}\alpha\quad\text{if}\quad\varepsilon_i=\half\\
 \alpha^{\overline{i}}\quad\text{if}\quad\varepsilon_i=0\end{array}\right..
 $$
Thus, the equation \eqref{theta} becomes
\begin{equation}\label{theta7}
\sum_{\varepsilon\in\mathcal{E},\alpha} D^\varepsilon_\alpha \Omega_\varepsilon^{\alpha}+\theta\sum_{\varepsilon\in\mathcal{O},\alpha} C^\varepsilon_\alpha \sum_{i=1}^n\Omega_\varepsilon^{\alpha^i}=-\sum_{\varepsilon\in\mathcal{O},\alpha} C^\varepsilon_\alpha\sum_{i=1}^n \xi^i_\varepsilon \Omega_{\varepsilon^i}^{\alpha_\varepsilon^i}-\theta\sum_{\varepsilon\in\mathcal{E},\alpha} D^\varepsilon_\alpha \sum_{i=1}^n\xi^i_\varepsilon \Omega_{\varepsilon^i}^{\alpha_\varepsilon^i}
\end{equation}
or equivalently (according to \eqref{varepsilon}, \eqref{ei} and \eqref{alpha})
\begin{equation}\label{theta8}\small{
\sum_{\varepsilon\in\mathcal{E},\alpha} D^\varepsilon_\alpha \Omega_\varepsilon^{\alpha}+\theta\sum_{\varepsilon\in\mathcal{O},\alpha}\left(\sum_{i=1}^n C^\varepsilon_{\alpha^{\overline{i}}} \right)\Omega_\varepsilon^{\alpha}=-\sum_{\varepsilon\in\mathcal{E},\alpha} \left(\sum_{i=1}^n\xi^i_\varepsilon C^{\varepsilon^i}_{\alpha^{\overline{i}}_\varepsilon}\right) \Omega_{\varepsilon}^{\alpha}-\theta\sum_{\varepsilon\in\mathcal{O},\alpha} \left(\sum_{i=1}^n\xi^i_\varepsilon D^{\varepsilon^i}_{\alpha^{\overline{i}}_\varepsilon}\right) \Omega_{\varepsilon}^{\alpha}}
\end{equation}
Thus, the coefficients $C^\varepsilon_\alpha$ must satisfy the following conditions
\begin{equation}\label{theta9}
\left\{\begin{array}{llll}D^\varepsilon_\alpha &=&-\sum_{i=1}^n\xi^i_\varepsilon C^{\varepsilon^i}_{\alpha^{\overline{i}}_\varepsilon}\quad\text{for}\quad \varepsilon\in\mathcal{E}\\
\sum_{i=1}^n C^\varepsilon_{\alpha^{\overline{i}}}&=&-\sum_{i=1}^n\xi^i_\varepsilon D^{\varepsilon^i}_{\alpha^{\overline{i}}_\varepsilon}\quad\text{for}\quad \varepsilon\in\mathcal{O}.\end{array}\right.
\end{equation}
Thus,
$$
\Omega(X_\theta)=\sum_{\varepsilon\in\mathcal{O},\alpha}C^{\varepsilon}_\alpha \Omega_{\varepsilon}^\alpha-
\theta\sum_{\varepsilon\in\mathcal{E},\alpha}\left(\sum_{i=1}^n\xi^i_\varepsilon C^{\varepsilon^i}_{\alpha^{\overline{i}}_\varepsilon}\right) \Omega_{\varepsilon}^\alpha=\sum_{\varepsilon\in\mathcal{O},\alpha}C^{\varepsilon}_\alpha \Omega_{\varepsilon}^\alpha-
\theta\sum_{\varepsilon\in\mathcal{O},\alpha} C^\varepsilon_\alpha\sum_{i=1}^n \xi^i_\varepsilon \Omega_{\varepsilon^i}^{\alpha_\varepsilon^i}
$$
with
\begin{equation}\label{alph}
\sum_{i=1}^n C^\varepsilon_{\alpha^{\overline{i}}}=\sum_{i,j}\xi^i_\varepsilon\xi^j_{\varepsilon^i} C^{\varepsilon^{ij}}_{\alpha^{\overline{ij}}_\varepsilon},
\end{equation}
where $\alpha^{\overline{ij}}_\varepsilon=(\alpha^{\overline{i}}_\varepsilon)_{\varepsilon^i}^{\overline{j}}$. But, it is easy to see that $$\xi^i_\varepsilon\xi^j_{\varepsilon^i}=-\xi^j_{\varepsilon}\xi^i_{\varepsilon^j},\quad \varepsilon^{ii}=\varepsilon\quad\text{and}\quad \alpha^{\overline{ii}}_\varepsilon=\alpha^{\overline{i}}.$$ Therefore the relation \eqref{alph} is trivial.
Of course, we agree that if  $\alpha^{\overline{i}}_\varepsilon\notin\mathbb{N}^n$ or  $\alpha^{\overline{i}}\notin\mathbb{N}^n$ then the correspondent coefficients $C^{\varepsilon^i}_{\alpha^{\overline{i}}_\varepsilon}$ or $C^{\varepsilon^i}_{\alpha^{\overline{i}}}$ are zero.

Now, let $B\in(\mathfrak{D}_{\lambda,\mu})_{\bar0}$, according to the decomposition \eqref{deven}, we can write $$B=\sum_{\varepsilon\in\mathcal{E},|\alpha|=\delta_\varepsilon+1}B^{\varepsilon}_\alpha \Omega_{\varepsilon}^\alpha+\theta\sum_{\varepsilon\in\mathcal{O},|\alpha|=\delta_\varepsilon+\frac{1}{2}}R^{\varepsilon}_\alpha \Omega_{\varepsilon}^\alpha.$$

\begin{equation}\label{theta10} \small{\partial B(X_\theta)=\displaystyle\half\sum_{\varepsilon\in\mathcal{O},|\alpha|=\delta_\varepsilon+\frac{1}{2}}\left(R^{\varepsilon}_\alpha-\sum_{i=1}^n\xi^i_\varepsilon B^{\varepsilon^i}_{\alpha^{\overline{i}}_\varepsilon}\right) \Omega_{\varepsilon}^\alpha+\half\theta\sum_{\varepsilon\in\mathcal{E},|\alpha|=\delta_\varepsilon+1}\left(\sum_{i=1}^n B^\varepsilon_{\alpha^{\overline{i}}} -\sum_{i=1}^n\xi^i_\varepsilon R^{\varepsilon^i}_{\alpha^{\overline{i}}_\varepsilon}\right) \Omega_{\varepsilon}^\alpha.}\end{equation}
Thus, according to Lemma \ref{lem2}, we have
$$\Omega=\partial B$$
where
$$B=2\theta\sum_{\varepsilon\in\mathcal{O},|\alpha|=\delta_\varepsilon+\frac{1}{2}}C^{\varepsilon}_\alpha \Omega_{\varepsilon}^\alpha.$$

Now, if $\Omega$ is an odd 1-cocycle then $\delta+\half\in\mathbb{N}$ and $\Omega(X_\theta)$ is even, so
$$
\Omega(X_\theta)=\sum_{\varepsilon\in\mathcal{E},|\alpha|=\delta_\varepsilon+\frac{1}{2}}C^{\varepsilon}_\alpha \Omega_{\varepsilon}^\alpha+
\theta\sum_{\varepsilon\in\mathcal{O},|\alpha|=\delta_\varepsilon+1}D^{\varepsilon}_\alpha \Omega_{\varepsilon}^\alpha.
$$
By the equation \eqref{sltr3} we have
\begin{equation}\label{theta11}\begin{array}{llll} D(\Omega(X_\theta,F_1,\dots,F_n))=\sum_{i=1}^n\xi_i
\Omega(X_\theta,F_1,\dots,D(F_i),\dots, F_n).\end{array}\end{equation}
Thus, the equation \eqref{theta} becomes
\begin{equation}\label{theta12}
\sum_{\varepsilon\in\mathcal{O},\alpha} D^\varepsilon_\alpha \Omega_\varepsilon^{\alpha}+\theta\sum_{\varepsilon\in\mathcal{E},\alpha} C^\varepsilon_\alpha \sum_{i=1}^n\Omega_\varepsilon^{\alpha^i}=\sum_{\varepsilon\in\mathcal{E},\alpha} C^\varepsilon_\alpha\sum_{i=1}^n \xi^i_\varepsilon \Omega_{\varepsilon^i}^{\alpha_\varepsilon^i}+\theta\sum_{\varepsilon\in\mathcal{O},\alpha} D^\varepsilon_\alpha \sum_{i=1}^n\xi^i_\varepsilon \Omega_{\varepsilon^i}^{\alpha_\varepsilon^i}
\end{equation}
or equivalently (according to \eqref{varepsilon} and\eqref{alpha})
\begin{equation}\label{theta13}
\sum_{\varepsilon\in\mathcal{O},\alpha} D^\varepsilon_\alpha \Omega_\varepsilon^{\alpha}+\theta\sum_{\varepsilon\in\mathcal{E},\alpha}\left(\sum_{i=1}^n C^\varepsilon_{\alpha^{\overline{i}}} \right)\Omega_\varepsilon^{\alpha}=\sum_{\varepsilon\in\mathcal{E},\alpha} \left(\sum_{i=1}^n\xi^i_\varepsilon C^{\varepsilon^i}_{\alpha^{\overline{i}}_\varepsilon}\right) \Omega_{\varepsilon}^{\alpha}+\theta\sum_{\varepsilon\in\mathcal{O},\alpha} \left(\sum_{i=1}^n\xi^i_\varepsilon D^{\varepsilon^i}_{\alpha^{\overline{i}}_\varepsilon}\right) \Omega_{\varepsilon}^{\alpha}
\end{equation}
Thus, the coefficients $C^\varepsilon_\alpha$ must satisfy the following conditions
\begin{equation}\label{theta14}
D^\varepsilon_\alpha =\sum_{i=1}^n\xi^i_\varepsilon C^{\varepsilon^i}_{\alpha^{\overline{i}}_\varepsilon}\quad\text{for}\quad \varepsilon\in\mathcal{O}.
\end{equation}
Thus,
$$
\Omega(X_\theta)=\sum_{\varepsilon\in\mathcal{E},\alpha}C^{\varepsilon}_\alpha \Omega_{\varepsilon}^\alpha+
\theta\sum_{\varepsilon\in\mathcal{O},\alpha}\left(\sum_{i=1}^n\xi^i_\varepsilon C^{\varepsilon^i}_{\alpha^{\overline{i}}_\varepsilon}\right) \Omega_{\varepsilon}^\alpha=\sum_{\varepsilon\in\mathcal{E},\alpha}C^{\varepsilon}_\alpha \Omega_{\varepsilon}^\alpha+
\theta\sum_{\varepsilon\in\mathcal{E},\alpha} C^\varepsilon_\alpha\sum_{i=1}^n \xi^i_\varepsilon \Omega_{\varepsilon^i}^{\alpha_\varepsilon^i}.
$$



Now, let $B\in(\mathfrak{D}_{\lambda,\mu})_{\bar1}$, according to the decomposition \eqref{deven}, we can write $$B=\sum_{\varepsilon\in\mathcal{O},|\alpha|=\delta_\varepsilon+1}B^{\varepsilon}_\alpha \Omega_{\varepsilon}^\alpha+\theta\sum_{\varepsilon\in\mathcal{E},|\alpha|=\delta_\varepsilon+\frac{1}{2}}R^{\varepsilon}_\alpha \Omega_{\varepsilon}^\alpha.$$

\begin{equation}\label{theta15} \partial B(X_\theta)=\displaystyle\half\sum_{\varepsilon\in\mathcal{O},\alpha}\left(R^{\varepsilon}_\alpha+\sum_{i=1}^n\xi^i_\varepsilon B^{\varepsilon^i}_{\alpha^{\overline{i}}_\varepsilon}\right) \Omega_{\varepsilon}^\alpha+\half\theta\sum_{\varepsilon\in\mathcal{E},\alpha}\left(\sum_{i=1}^n B^\varepsilon_{\alpha^{\overline{i}}} +\sum_{i=1}^n\xi^i_\varepsilon R^{\varepsilon^i}_{\alpha^{\overline{i}}_\varepsilon}\right) \Omega_{\varepsilon}^\alpha.\end{equation}
Thus, according to Lemma \ref{lem2}, we have
$$\Omega=\partial B$$
where
$$B=2\theta\sum_{\varepsilon\in\mathcal{E},|\alpha|=\delta_\varepsilon+\frac{1}{2}}C^{\varepsilon}_\alpha \Omega_{\varepsilon}^\alpha.$$

\hfill$\Box$
\end{proofname}
\begin{corollary}\label{sa}
Any 1-cocycle $\Omega\in Z^1_\mathrm{diff}(\frak{aff}(1|1),\frak{D}_{\lambda,\mu})$ of $\mathfrak{aff}(1|1)$ is a coboundary if
and only if its restriction  to $\mathfrak{aff}(1)$ is
a coboundary.
\end{corollary}
 \section{The space $\HH^1_\mathrm{diff}(\frak
{aff}(1|1),\frak{D}_{\l,\mu})$}
According to Theorem \ref{cor} the restriction of any 1-cocycle $\Omega\in Z^1_\mathrm{diff}(\frak{aff}(1|1),\frak{D}_{\lambda,\mu})$ to $\mathfrak{aff}(1)$  has the following structure:
\begin{proposition}\label{prop}\begin{enumerate}
                     \item [1)] If $2\delta\notin\mathbb{N}$ then $\HH^1_\mathrm{diff}(\frak
{aff}(1|1),\frak{D}_{\l,\mu})=0$.
                     \item [2)] If $\delta\in\mathbb{N}$ then, up to a coboundary,  any  1-cocycle $\Omega\in Z^1_\mathrm{diff}(\frak{aff}(1|1),\frak{D}_{\lambda,\mu})$ is even and its restriction to $\mathfrak{aff}(1)$  has the following form:
$$
\Omega(X_h)=\sum_{\varepsilon\in\mathcal{E},|\alpha|=\delta_\varepsilon}C^{\varepsilon}_\alpha h' \Omega_{\varepsilon}^\alpha+
\theta\sum_{\varepsilon\in\mathcal{O},|\alpha|=\delta_\varepsilon+\half}D^{\varepsilon}_\alpha h' \Omega_{\varepsilon}^\alpha
$$
where $\delta_\varepsilon=\mu-|\lambda+\varepsilon|$.
                     \item [3)] If $\delta\in\mathbb{N}+\half$ then, up to a coboundary, any  1-cocycle $\Omega\in Z^1_\mathrm{diff}(\frak{aff}(1|1),\frak{D}_{\lambda,\mu})$ is odd and its restriction to $\mathfrak{aff}(1)$ has the following form:
$$
\Omega(X_h)=\sum_{\varepsilon\in\mathcal{O},|\alpha|=\delta_\varepsilon}C^{\varepsilon}_\alpha h' \Omega_{\varepsilon}^\alpha+
\theta\sum_{\varepsilon\in\mathcal{E},|\alpha|=\delta_\varepsilon+\half}D^{\varepsilon}_\alpha h' \Omega_{\varepsilon}^\alpha.
$$
                   \end{enumerate}

\end{proposition}
\begin{thm}\begin{enumerate}
                     \item [1)] If $2\delta\notin\mathbb{N}$ then $\HH^1_\mathrm{diff}(\frak
{aff}(1|1),\frak{D}_{\l,\mu})=0$.
\item [2)] If $\delta=k\in\mathbb{N}$ then
 $$\mathrm{dim}\mathrm{H}^1(\mathfrak{aff}(1|1);\mathfrak{D}_{\lambda,\mu})=\sum_{r=0}^{[\frac{n}{2}]}\begin{pmatrix}n\\2r\end{pmatrix} \begin{pmatrix}n+k-r-1\\k-r\end{pmatrix}.$$
  \item [3)] If $\delta=k+\frac{1}{2}\in\mathbb{N}+\frac{1}{2}$ then
 $$\mathrm{dim}\mathrm{H}^1(\mathfrak{aff}(1|1);\mathfrak{D}_{\lambda,\mu})=\sum_{r=0}^{[\frac{n-1}{2}]}\begin{pmatrix}n\\2r+1\end{pmatrix} \begin{pmatrix}n+k-r-1\\k-r\end{pmatrix}.$$
 \end{enumerate}
 A basis of $\mathrm{H}^1(\mathfrak{aff}(1|1);\mathfrak{D}_{\lambda,\mu})$ is given by the family $\Gamma^\alpha_\varepsilon$ defined on $\mathfrak{aff}(1)$ by $$\Gamma^\alpha_\varepsilon(X_h)=h'\Omega^\alpha_\varepsilon+\theta\sum_{i=1}^n\xi^i_\varepsilon h' \Omega_{\varepsilon^i}^{\alpha^i_\varepsilon}$$ where $\varepsilon\in\mathcal{E}$ and $|\alpha|=\delta_\varepsilon$ if $\delta\in\mathbb{N}$ and $\varepsilon\in\mathcal{O}$ and $|\alpha|=\delta_\varepsilon+\frac{1}{2}$ if $\delta\in\mathbb{N}+\frac{1}{2}$.
\end{thm}

\begin{proofname}. If $\delta\in\mathbb{N}$ then, according to Proposition \ref{prop}, we have for any $X_h\in\mathfrak{aff}(1)$
$$
\Omega(X_h)=\sum_{\varepsilon\in\mathcal{E},|\alpha|=\delta_\varepsilon}C^{\varepsilon}_\alpha h' \Omega_{\varepsilon}^\alpha+
\theta\sum_{\varepsilon\in\mathcal{O},|\alpha|=\delta_\varepsilon+\half}D^{\varepsilon}_\alpha h' \Omega_{\varepsilon}^\alpha.
$$
In this case, $X_\theta$ is an odd operator, so, it is of the following form:
$$
\Omega(X_\theta)=\sum_{\varepsilon\in\mathcal{O},\alpha}B^{\varepsilon}_\alpha  \Omega_{\varepsilon}^\alpha+
\theta\sum_{\varepsilon\in\mathcal{E},\alpha}R^{\varepsilon}_\alpha  \Omega_{\varepsilon}^\alpha.
$$
Using the 1-cocycle conditions we seek to establish relationships that must be satisfied by the coefficients $B^{\varepsilon}_\alpha $, $R^{\varepsilon}_\alpha $, $C^{\varepsilon}_\alpha $, \dots

The 1-cocycle condition reads:
\begin{equation}
\begin{array}{lll}\label{part1}
X_x\cdot\Omega(X_\theta)-{X_\theta}\cdot\Omega(X_x)+\half
\Omega(X_\theta)=0,
\end{array}
\end{equation}
where
$$\begin{array}{lllllllll}
X_x\cdot\Omega(X_\theta)&=&\mathfrak{L}^\mu_{X_x}\circ\Omega(X_\theta)-\Omega(X_\theta)\circ\mathfrak{L}^\lambda_{X_x},\\
\Omega_{\varepsilon}^\alpha\circ\mathfrak{L}^\lambda_{X_x}&=&\sum_{i=1}^n\left(x\Omega_{\varepsilon}^{\alpha^i}+(\alpha_i+\lambda_i+\varepsilon_i) \Omega_{\varepsilon}^\alpha\right),\\
\mathfrak{L}^\mu_{X_x}\circ\Omega(X_\theta)&=&\sum_{\varepsilon\in\mathcal{O},\alpha}B^{\varepsilon}_\alpha L^\mu_{X_x} \Omega_{\varepsilon}^\alpha+
\theta\sum_{\varepsilon\in\mathcal{E},\alpha}R^{\varepsilon}_\alpha L^{\mu+\half}_{X_x} \Omega_{\varepsilon}^\alpha,\\
L^\mu_{X_x} \Omega_{\varepsilon}^\alpha&=&x\sum_{i=1}^n \Omega_{\varepsilon}^{\alpha^i}+\mu \Omega_{\varepsilon}^\alpha,\\{2{X_\theta}\cdot\Omega(X_x)}&=&{\sum_{\varepsilon\in\mathcal{O},\alpha} D^\varepsilon_\alpha \Omega_\varepsilon^{\alpha}+\theta\sum_{\varepsilon\in\mathcal{E},\alpha} C^\varepsilon_\alpha \sum_{i=1}^n\Omega_\varepsilon^{\alpha^i}}-\\&~&\sum_{\varepsilon\in\mathcal{E},\alpha} C^\varepsilon_\alpha\sum_{i=1}^n \xi^i_\varepsilon \Omega_{\varepsilon^i}^{\alpha_\varepsilon^i}-\theta\sum_{\varepsilon\in\mathcal{O},\alpha} D^\varepsilon_\alpha \sum_{i=1}^n\xi^i_\varepsilon \Omega_{\varepsilon^i}^{\alpha_\varepsilon^i}\end{array}
$$
Therefore, the relation \eqref{part1} is satisfied if and only if
$$\left\{\begin{array}{ll}
(\delta_\varepsilon-|\alpha|+\half) B^{\varepsilon}_\alpha-\half D^\varepsilon_\alpha+\half\left(\sum_{i=1}^n\xi^i_\varepsilon C^{\varepsilon^i}_{\alpha^{\overline{i}}_\varepsilon}\right)=0\\
(\delta_\varepsilon-|\alpha|+1) R^{\varepsilon}_\alpha-\half\sum_{i=1}^n C^{\varepsilon}_{\alpha^{\overline{i}}}+\half\sum_{i=1}^n\xi^i_\varepsilon D^{\varepsilon^i}_{\alpha^{\overline{i}}_\varepsilon}=0.\end{array}\right.
$$
Thus, for $|\alpha|=\delta_\varepsilon+\half$ we have
$$
D^\varepsilon_\alpha=\sum_{i=1}^n\xi^i_\varepsilon C^{\varepsilon^i}_{\alpha^{\overline{i}}_\varepsilon}
$$
and for $|\alpha|\neq\delta_\varepsilon+\half$ we have $|\alpha^{\overline{i}}_\varepsilon|\neq\delta_{\varepsilon^i}$, therefore $D^\varepsilon_\alpha=C^{\varepsilon^i}_{\alpha^{\overline{i}}_\varepsilon}=0$ (indeed $D^\varepsilon_\alpha$ appear only for $|\alpha|=\delta_\varepsilon+\half$ and $C^{\varepsilon}_{\alpha}$ appear only for $|\alpha|=\delta_\varepsilon$).
Therefore, $B^{\varepsilon}_\alpha=0$ for $|\alpha|\neq\delta_\varepsilon+\half$.

For $|\alpha|=\delta_\varepsilon+1$ we have
$$
\sum_{i=1}^n C^{\varepsilon}_{\alpha^{\overline{i}}}=\sum_{i=1}^n\xi^i_\varepsilon D^{\varepsilon^i}_{\alpha^{\overline{i}}_\varepsilon}
$$
and for $|\alpha|\neq\delta_\varepsilon+1$ we have $C^{\varepsilon}_{\alpha^{\overline{i}}}=D^{\varepsilon^i}_{\alpha^{\overline{i}}_\varepsilon}=0$, therefore $R^{\varepsilon}_\alpha=0$.

Thus,
$$
\Omega(X_\theta)=\sum_{\varepsilon\in\mathcal{O},|\alpha|=\delta_\varepsilon+\half}B^{\varepsilon}_\alpha  \Omega_{\varepsilon}^\alpha+
\theta\sum_{\varepsilon\in\mathcal{E},|\alpha|=\delta_\varepsilon+1}R^{\varepsilon}_\alpha  \Omega_{\varepsilon}^\alpha
$$
and
$$
\Omega(X_h)=\sum_{\varepsilon\in\mathcal{E},|\alpha|=\delta_\varepsilon}C^{\varepsilon}_\alpha h' \Omega_{\varepsilon}^\alpha+
\theta\sum_{\varepsilon\in\mathcal{O},|\alpha|=\delta_\varepsilon+\half}\sum_{i=1}^n\xi^i_\varepsilon C^{\varepsilon^i}_{\alpha^{\overline{i}}_\varepsilon} h' \Omega_{\varepsilon}^\alpha
$$
From Corollary \ref{sa} that the dimension of $\mathrm{H}^1(\mathfrak{aff}(1|1);\mathfrak{D}_{\lambda,\mu})$ is equal to the number of parameters $C^{\varepsilon}_\alpha$. That is, $$\mathrm{dim}\mathrm{H}^1(\mathfrak{aff}(1|1);\mathrm{D}_{\lambda,\mu})=\sum_{r=0}^{[\frac{n}{2}]}\begin{pmatrix}n\\2r\end{pmatrix} \begin{pmatrix}n+k-r-1\\k-r\end{pmatrix},$$ where $\begin{pmatrix}n\\2r\end{pmatrix}$  is the number of $\varepsilon$ with $|\varepsilon|=r$ and $\begin{pmatrix}n+k-r-1\\k-r\end{pmatrix}$ is the number of $\alpha$ with $|\alpha|=\delta_\varepsilon=k-r$.

\bigskip

Now, If $\delta\in\mathbb{N}+\frac{1}{2}$ then  we have
$$
\Omega(X_h)=\sum_{\varepsilon\in\mathcal{O},|\alpha|=\delta_\varepsilon}C^{\varepsilon}_\alpha h' \Omega_{\varepsilon}^\alpha+
\theta\sum_{\varepsilon\in\mathcal{E},|\alpha|=\delta_\varepsilon+\half}D^{\varepsilon}_\alpha h' \Omega_{\varepsilon}^\alpha.
$$
In this case, $X_\theta$ is an even operator, so, it is of the following form:
$$
\Omega(X_\theta)=\sum_{\varepsilon\in\mathcal{E},\alpha}B^{\varepsilon}_\alpha  \Omega_{\varepsilon}^\alpha+
\theta\sum_{\varepsilon\in\mathcal{O},\alpha}R^{\varepsilon}_\alpha  \Omega_{\varepsilon}^\alpha.
$$
As before we prove that
$$\left\{\begin{array}{ll}
(\delta_\varepsilon-|\alpha|+\half) B^{\varepsilon}_\alpha-\half D^\varepsilon_\alpha+\half\left(\sum_{i=1}^n\xi^i_\varepsilon C^{\varepsilon^i}_{\alpha^{\overline{i}}_\varepsilon}\right)=0\\
(\delta_\varepsilon-|\alpha|+1) R^{\varepsilon}_\alpha-\half\sum_{i=1}^n C^{\varepsilon}_{\alpha^{\overline{i}}}+\half\sum_{i=1}^n\xi^i_\varepsilon D^{\varepsilon^i}_{\alpha^{\overline{i}}_\varepsilon}=0.\end{array}\right.
$$
Thus, $$\mathrm{dim}\mathrm{H}^1(\mathfrak{aff}(1|1);\mathfrak{D}_{\lambda,\mu})=\sum_{r=0}^{[\frac{n-1}{2}]}\begin{pmatrix}n\\2r+1\end{pmatrix} \begin{pmatrix}n+k-r-1\\k-r\end{pmatrix},$$ where  $\begin{pmatrix}n\\2r+1\end{pmatrix}$ is the number of $\varepsilon$ with  $|\varepsilon|=r+\frac{1}{2}$ and $\begin{pmatrix}n+k-r-1\\k-r\end{pmatrix}$ is the number of $\alpha$ with $|\alpha|=k-r$.

\hfill$\Box$
\end{proofname}

For $n=2$, $\mathrm{dim}\mathrm{H}^1(\mathfrak{aff}(1|1);\mathfrak{D}_{\lambda,\mu})=2\delta+1$.


\begin{thebibliography}{99}

\bibitem{ab}
Agrebaoui B and Ben Fraj N, On the cohomology of the Lie
superalgebra of contact vector fields on $S^{1|1}$, {\it Bell.
Soc. Roy. Sci. Li\`ege} {\bf 72}, 6, 2004, 365--375.

\bibitem{abs}
Arnal D, Ben Ammar M and Selmi M, Normalisation d'une repr\'esentation non lin\'eaire d'une alg\`ebre de Lie, {\it Annales de la facult\'e des sciences de Toulouse, 5$^e$ s\'erie, tome 9, n$^o$3, (1988), p 355--579.}

\bibitem{bab}
Basdouri I, Ben Ammar M, Cohomology of $\mathfrak {osp}(1|2)$
acting on linear differential operators on the supercircle
$S^{1|1}$. Letters in Mathematical Physics (2007) 81:239--251.

 \bibitem{bjs}
Ben Ammar M, Jabeur A, Safi I, cohomology of $\frak {osp}(1|2)$ acting
on the space of bilinear differential operators on the superspace
$\mathbb{R}^{1|1}$, International Journal of Geometric Methods in Modern Physics
 (2000).
\bibitem{b}
Bouarroudj S, {\it Cohomology of the vector fields Lie algebras
on $\mathbb{R}\mathbb{P}^1$ acting on bilinear differential
operators}, International Journal of Geometric Methods in Modern Physics
 (2005), {\bf 2}; N 1,  23-40.
\bibitem{Fu}
Fuchs  D B, {\it Cohomology of infinite-dimensional Lie algebras},
Plenum Publ. New York, 1986.
\bibitem{g}
Gargoubi H, {\it Sur la g\'eom\'etrie de l'espace des op\'erateurs
diff\'erentiels lin\'eaires sur $\mathbb{R}$}, Bull. Soc. Roy. Sci.
Li\`ege. Vol. 69, 1, 2000, 21�47.
\bibitem{gmo} Gargoubi H, Mellouli N
and Ovsienko V {\it Differential operators on supercircle:
conformally equivariant quantization and symbol calculus,} Letters
in Mathematical Physics (2007) {\bf 79}:51�65.
\bibitem{lec}Lecomte B P A, {\it On the cohomology of
$\frak{sl}(n + 1;\mathbb{R})$ acting on differential operators and
$\frak{sl}(n + 1;\mathbb{R})$-equivariant symbols,} Indag. Math. NS.
11 (1), (2000), 95 114.


\bibitem{r1}
A. Nijenhuis,  R.W. Richardson, {\it Cohomology and deformations
in graded Lie algebras}. Bull. Amer. Math. Soc., 72 (1966), 1-29.

\end{thebibliography}
\end{document}